\documentclass {article}
\newcommand {\ncm }{\newcommand }
\ncm {\bma }  {\begin  {displaymath}}
\ncm {\ema }  {\end    {displaymath}}
\ncm {\beq }  {\begin  {equation}}
\ncm {\eeq }  {\end    {equation}}
\ncm {\bea }  {\begin  {eqnarray}}
\ncm {\eea }  {\end    {eqnarray}}
\ncm {\beast }{\begin  {eqnarray*}}
\ncm {\eeast }{\end    {eqnarray*}}

\title {The~Fundamental~Theorem~of~Geometric~Calculus via a Generalized Riemann Integral}
\author
{Alan Macdonald  \\
Department of Mathematics \\
Luther College, Decorah, IA 52101,  U.S.A.\\ 
macdonal@luther.edu}

\begin {document}
\maketitle

\begin{abstract}
\noindent
Using recent advances in integration theory, we give a proof of the 
fundamental theorem of geometric calculus. 
We assume only that the tangential derivative $\nabla \! _VF$ exists 
and is Lebesgue integrable. We also give sufficient conditions that $\nabla \! _VF$ exists. 
\end{abstract}

\vspace{.15in}

\vspace{.02in}
\normalsize
\vspace{.2in}

{\bf I. Introduction.} We assume that the reader is familiar with the fundamentals of geometric calculus \cite{gc1,gc2}.  We prove here this version of the

\vspace{.1in}
\noindent
\textbf { Fundamental Theorem of Geometric Calculus.}
\emph{
Let $M $ be a $C^1$ $n $--chain in ${\mathbf R}^m $. 
Suppose that the multivector field $F $ is continuous on $M $ and that the tangential derivative  $\nabla \! _VF$ (defined below) exists and is Lebesgue integrable on $M-\partial M$. Then }
\beq
\int_M {dV\; \nabla \! _VF}=\int_{\partial M} {dA\,F.}  \label {eq:ftcgc}
\eeq

Here $V $ is the tangent to $M$ and $A $ is the tangent to $\partial M $. (By \emph{the} tangent, we mean, e.g., that $V(X)$ is the unit positively oriented pseudoscalar in the tangent algebra to $M$ at $X$.) Recall the important relationship $V{^\dagger\! }A = N $, where $N $ is the unit outward normal to $M $ \cite[p. 319]{gc1}.  The relationships $dV = |dV|V $ and $dA = |dA|A $ define the integrals componentwise as Lebesgue integrals on $M $ and $\partial M $ \cite[p. 317]{gc1}. 

Eq. \ref {eq:ftcgc} is a generalization of the fundamental theorem of calculus and the integral theorems of vector calculus \cite[p. 323]{gc1}, Cauchy's theorem \cite{mf}, and an arbitrary dimension multivector version of Cauchy's theorem \cite{bds,mf}. 

\newpage
Hestenes gives a heuristic demonstration of the fundamental theorem \cite{gc1}.  The demonstration is wonderfully brief and offers great intuitive insight, but it is not a rigorous proof. In particular, the demonstration does not provide sufficient conditions on $F$ for Eq. \ref {eq:ftcgc} to hold. Some condition is necessary even for the fundamental theorem of calculus. A standard example is the function 
$f(x)=x^2\cos (\pi /x^2)$ for $x \in (0,1]$, with $f(0) = 0$. The derivative $f'$ exists on $[\,0, 1]$. But $f$ is not absolutely continuous on $[\,0, 1]$, and so $f'$ is not Lebesgue integrable there. We show that Eq. \ref {eq:ftcgc} holds whenever $\nabla \!_VF$ is Lebesgue integrable.

This may come as a surprise, as most similar theorems assume more than Lebesgue integrability of a derivative. Some assume that the derivative is continuous. Examples are most versions of Stokes' theorem on manifolds, Cauchy's theorem, and the abovementioned multivector version of Cauchy's theorem. Others assume that the derivative is Riemann integrable. Examples are the one dimensional fundamental theorem of calculus and a recent version of Stokes' theorem \cite{acker}. 

We shall prove the multivector version of Cauchy's theorem as a corollary to our fundamental theorem assuming only that the derivative exists. This is a generalization of the Cauchy-Goursat theorem, not just the Cauchy theorem.

We shall prove the fundamental theorem first for $M= [\,0,\;1]^n$, and then for $M $ a $C^1$ $n$--chain. (This is against the spirit of Hestenes' appealing program to work in vector manifolds in a coordinate free manner \cite{dfgc,gc2}. Perhaps I suffer coordinitis \cite{mathvirus}, but I see no way to obtain the results of this paper within Hestenes' program.)

Let $f $ be a multivector valued function defined in a neighborhood $u $ of $x\in {\mathbf R}^n$. We define the {\em gradient}, $\nabla \! f(x)$. (This is the tangential derivative
$\nabla \! _vf(x)$ in ${\mathbf R}^n $ where $v(x)\equiv v$, the unit positively oriented pseudoscalar in ${\mathbf R}^n$.) Let $c $ be a cube with $x\in c\subseteq u$. Define
\beq
\nabla \! f\!\left( x \right)\;\;=\;\;v^\dagger\!  \lim \limits_{\scriptstyle {x\in c}\atop
  \scriptstyle {{\rm diam} \left( {c } \right)\to 0}}{1 \over {\left| {\,c\,} \right|}}\int_{\partial c} {da\,f}.    \label {eq:grad}
\eeq

If $f $ is differentiable, i.e., if each component of $f $ is differentiable, 
then $\nabla \!  f $ exists and $\nabla \!  f = \sum\nolimits_i {e_i\partial _if}$ in every rectangular coordinate system. This is a special case of the theorem of Sec. V.

The definition of $\nabla \!  f $, Eq. \ref {eq:grad}, is the key to Hestenes' heuristic demonstration of the fundamental theorem. It is also a key to our proof of the theorem. There is a second key to our proof, again a definition, this time of an integral.

\vspace{.15in}
{\bf II. The RP integral.} The RP integral \cite{Mawhin1, Mawhin2} is one of several {\em generalized Riemann integrals} which were designed to overcome the deficiency of the Lebesgue integral shown by the example above: the Lebesgue integral does not integrate all derivatives. In ${\mathbf R}^1$ this deficiency is completely removed by the {\em Henstock-Kurzweil} (HK) generalized Riemann integral: 
If $f' $ exists on an interval $[a,b]$, then $f'$ is HK integrable over $[a,b]$ to $f(b)-f(a)$ \cite{Bartle, McLeod}. The HK integral is {\em super Lebesgue}: If $f$ is Lebesgue integrable, then it is HK integrable to the same value.

The HK integral in ${\mathbf R}^n$ is super Lebesgue, but it does not integrate all divergences \cite[Example 5.7] {Pfeffer1}.  The RP integral was designed to overcome this \cite[Theorem 2]{Mawhin1}. It is also super Lebesgue. To see this, note that the definition of the HK integral on $[\,0, 1]^n$ is the same as that of the RP integral, except that $[\,0, 1]^n$ is partitioned into rectangles rather than cubes \cite[p. 33]{McLeod}. It follows that the RP integral is super HK. Also, the HK integral is super Lebesgue \cite[p. 236]{McLeod}. Thus the RP integral is super Lebesgue.

We shall prove that if $\nabla \! f $ exists on $[\,0,1]^n$, then it is RP integrable there.

If the RP integral is super Lebesgue and always integrates $\nabla \!  f $, why don't we abandon the Lebesgue integral in favor of the RP integral? Most important for us, the change of variable theorem fails \cite[p. 143]{Pfeffer2}, and so the integral cannot be lifted to manifolds. Fubini's theorem also fails \cite[Remark 5.8]{Pfeffer1}. And there are other deficiencies \cite[Remark 7.3]{Pfeffer1}.
Unlike ${\mathbf R}^1$, where the HK integral seems to be completely satisfactory, none of the several generalized Riemann integrals in higher dimensions has enough desirable properties to make it a useful general purpose integral. We use the RP integral here only as a catalyst: since the RP integral is super Lebesgue, we use it to compute Lebesgue integrals. This will allow us to make Hestenes' heuristic demonstration of the fundamental theorem rigorous.
 
\vspace{.1in}
We now give a series of definitions leading to the RP integral, specialized to $[\,0,1]^n$.

A {\em gauge} on $[\,0,1]^n$ is a positive function $\delta (x)$ on $[\,0,1]^n$. 
(The gauge function $\delta (x)$ in generalized Riemann integrals replaces the constant norm $\delta $ of a partition in the Riemann integral.)

A {\em tagged RP partition} $\{ c_j,\ x_j\} _{j=1}^k$ of $[\,0,1]^n$ is a decomposition of $[\,0,1]^n$ into closed subcubes $c_j$ together with points $x_j\in c_j$. The $c_j$ are disjoint except for boundaries and have sides parallel to the axes.

Let $\delta$ be a gauge on $[\,0,1]^n$. A tagged RP partition $\{ c_j,\ x_j\} _{j=1}^k$ is {\em $\delta$-fine} if ${\rm diam}(c_j)\le \delta (x_j),j=1\ldots k$.

Let $f$ be a multivector valued function defined on $[\,0,1]^n$. A multivector, denoted (RP)$\int_{\,\left[ {0,1} \right]^{n}} {dv\,f}$, is the {\em directed RP integral of \textrm{f} over $[\,0,1]^n$} if, given $\epsilon > 0$, there is a gauge $\delta$ on 
$[\,0,1]^n$ such that for every $\delta$-fine tagged RP partition of $[\,0,1]^n$,
\beq
\left| {\;\left( \rm {RP} \right)\int_{\,\left[ {0,1} \right]^n} {dv\,f}-\;v\sum\limits_{j=1}^k {|\,c_j|\,f(x_j)}\;} \right|\le \varepsilon .  \label {eq:RP}
\eeq

The definition makes sense only if a $\delta$-fine RP partition of $[\,0,1]^n$ exists. We now prove this (Cousin's lemma). First note that if a cube is partitioned into subcubes, each of which has a $\delta$-fine RP partition, then the original cube also has a $\delta$-fine RP partition. Thus if there is no $\delta$-fine RP partition of $[\,0,1]^n$, then there is a sequence $c_1 \supseteq c_2\supseteq \ldots $ of compact subcubes of $[\,0,1]^n$, with no $\delta$-fine RP partition and with ${\rm diam}\;c_i\to 0$. Let 
$\{ x\} =\bigcap _ic_i$. 
Choose $j$ so large that ${\rm diam}\;(c_j)\le \delta (x)$. Then $\{ (x,c_j)\} $
is a $\delta$-fine RP partition of $c_j$, which is a contradiction. (It is interesting to note that the proof of the Cauchy-Goursat theorem uses a similar compactness argument to obtain a contradiction.)

If the RP integral exists, then it is unique. This follows from the fact that if 
$\delta _1$ and $\delta _2$ are gauges and $\delta = \rm{Min}(\delta _1, \delta _2)$, then a $\delta$-fine RP partition is also $\delta _1$-fine and $\delta _2$-fine.

{\bf III. Proof of the Fundamental Theorem on} \mbox{\boldmath $[\,0,1]^n\!.$ }
On $[\,0,1]^n\!$ the fundamental theorem, Eq. \ref {eq:ftcgc}, becomes
\beq
\int_{\left[ {0,1} \right]^n} {dv\,\nabla \!  f}=\int_{\partial \left[ {0,1} \right]^n} {da\,f,}	 	\label {eq:e6}
\eeq
where $f$ is continuous on $[\,0,1]^n$ and $\nabla f$ exists and is Lebesgue integrable on $(0,1)^n$.

First suppose that $\nabla \!  f$ exists on all of $[\,0,1]^n$. (Thus $f$ is defined on an open set containing $[\,0,1]^n$.) Make no assumption about the integrability of $\nabla \! f$. We show that $\nabla \! f$ is RP-integrable.
Let $\varepsilon >0$ be given. Define a gauge $\delta$ on $[\,0,1]^n$ as follows. Choose $x \in [\,0,1]^n$. By the definition of $\nabla \!  f(x)$, Eq. \ref {eq:grad}, there is a $\delta (x) > 0$ so that if $c$ is a cube with $x \in c$ and ${\rm diam}(c) \le \delta (x)$ then
\bma
\left| {\,{{v^\dagger\!  } \over {\left| {\,c\,} \right|}}\int_{\partial \,c} {da\,f}-\nabla \!  f\left( x \right)\,} \right|\le \varepsilon \,.  
\ema
Let $\{ c_j,\ x_j\} _{j=1}^k$ be a $\delta$-fine RP partition of $[\,0,1]^n$. Then
\beast
\lefteqn{ \left| {\int_{\partial \left[ {0,1} \right]^n} {daf}-v\sum\limits_j {\left| {c_j} \right|\nabla \!  f\left( {x_j} \right)}} \right| } \\
& = & \left| {\sum\limits_j {\int_{\partial c_j} {daf}}-v\sum\limits_j {\left| {c_j} \right|\nabla \!  f(x_j)}} \right|  \\
&\le& \sum\limits_j {\left| {\int_{\partial c_j} {daf}-v\left| {c_j} \right|\nabla \!  f(x_j)} \right|}  \\
& = & \sum\limits_j {\left| {{{v^\dagger } \over {\left| {c_j} \right|}}\int_{\partial c_j} {daf}-\nabla \!  f(x_j)} \right|}\left| {c_j} \right|  \\
&\le & \sum\limits_j \varepsilon \left| {c_j} \right|=\varepsilon .
\eeast
Thus by the definition of the RP integral, Eq. \ref {eq:RP},
\bea
\left( \rm{RP} \right)\int_{\left[ {0,1} \right]^n} {dv\,\nabla \!  f}=\int_{\partial \left[ {0,1} \right]^n} {da\,f.}  \label {eq:e4}
\eea
The definitions of the gradient and the RP integral fit hand in glove to produce this equation. The equation shows that the RP integral integrates all gradients on $[\,0,1]^n$, and integrates them to the "right" value. 

Now suppose that $f$ is continuous on $[\,0,1]^n$, and that $\nabla \!  f$ exists and is Lebesgue integrable on $(0,1)^n$. Let $c_j = [j^{-1},1-j^{-1}]^n$. Then by Eq. \ref {eq:e4} (which applies to $c_j$),
\beq
\int_{c_j} {dv\,\nabla \!  \,f}=\int_{\partial c_j} {da\,f.}
                                                            \label {eq:e5}
\eeq
We have removed the RP designation since $\nabla \!  f$ is Lebesgue integrable on $c_j$.  

Let $j\to \infty $ in Eq. \ref {eq:e5}. By the Lebesgue dominated convergence theorem (applied componentwise), the left side of Eq. \ref {eq:e5} approaches the left side of Eq. \ref {eq:e6}. And by the uniform continuity of $f$ on $[\,0,1]^n$, the right side of Eq. \ref {eq:e5} approaches the right side of Eq. \ref {eq:e6}. This completes the proof.
\vspace{.1 in}

{\bf IV. Proof of the Fundamental Theorem on a} \mbox{\boldmath $C^1\ n$}{\bf --chain} \mbox{\boldmath $M$.} For convenience, we reproduce the statement of the theorem, Eq. \ref {eq:ftcgc}:
\bma
\int_M {dV\; \nabla \! _VF}=\int_{\partial M} {dA\,F.} 
\ema

We must first define $\nabla \!  _VF(X)$ in this equation, where $F$ is an ${\mathbf R}^m$-multivector valued function defined in a neighborhood $U$ of $X \in M \subseteq {\mathbf R}^m$. Let $c \subseteq [\,0,1]^n$ be a cube with $X \in \varphi(c) \equiv C \subseteq U$. Then 
\bma
\nabla \!  _VF\left( X \right)=V^\dagger\!  \left( X \right)\lim \limits_{\scriptstyle {X\in C}\atop
  \scriptstyle {{\rm diam}\left( c \right)\to 0}}{1 \over {\left| C \right|}}\int_{\partial C} {dA\,F}.
\ema
(Cf. the definition of the gradient, Eq. \ref {eq:grad}.) With $X = \varphi (x)$,
\bea
\nabla \!  _VF\left( X \right) & = & V^\dagger\!  \left( X \right)\lim \limits_{\scriptstyle {x\in c}\atop
  \scriptstyle {{\rm diam}\left( {c} \right)\to 0}}{{\left| {\,c\,} \right|} \over {\left| {\,C\,} \right|}}\;{1 \over {\left| {\,c\,} \right|}}\int_{\partial C} {dA\,F}  \nonumber \\
& = & {{V^\dagger\!  \left( X \right)} \over {J\left( x \right)}}\; \lim \limits_{\scriptstyle {x\in c}\atop
  \scriptstyle {{\rm diam} \left( {c} \right)\to 0}}{1 \over {\left| {\,c\,} \right|}}\int_{\partial C} {dA\,F}.  \label {eq:e7}
\eea

In Sec. V we show that if $F$ is differentiable on a $C^2\ n$-chain, then $\nabla \!  _VF$ exists and $\nabla \!  _VF=\sum {N_i\partial _iF}$ for certain vectors $N_i$.
\vspace{.1in}

{\bf Proof.} It is sufficient to prove the theorem for a singular $n$-cube $\varphi \!: [\,0,1]^n \rightarrow M \subseteq {\mathbf R}^m$, where $\varphi$ is continuous on $[\,0,1]^n$, and $\varphi '$, the differential of $\varphi$, exists and is continuous on $(0,1)^n$. As in Sec. III, we can further reduce our theorem to the case where $\nabla \!  _VF$ exists on all of $M$. In this reduced case $\varphi '$ is exists and is continuous on $[\,0, 1]^n$. Thus, with $J = {\rm det}(\varphi ') > 0$, ${\rm sup} J(x) < \infty$.

The Lebesgue integral $\int_M {dV\nabla \!  _VF}$ can be expressed, using $dV = \left| dV \right| V$ and the change of variable formula, as
\beq
\int_M {d V\,\nabla \!  _VF(X)}={\rm (RP)}\int_{\,\left[ {0,1} \right]^n} {\left| {\,dv\,} \right|\,J(x)\; V(X) \,\nabla \!  _VF(X).} \label   {eq:e9}
\eeq
The RP designation is allowed in the integral on the right side since the RP integral is super Lebesgue. We now compute this integral.

Let $\epsilon > 0$ be given. We define a gauge $\delta $ on $[\,0,1]^n$. Choose $x \in [\,0,1]^n$. From Eq. \ref {eq:e7}, there is a $\delta (x) > 0$ so that if $c$ is a cube with $x \in c$ and ${\rm diam}(c) \le \delta (x)$ then
\beq
\left| {\,{{V^\dagger\!  \left( X \right)} \over {J\left( x \right)\;\left| {\,c\,} \right|}}\;\int_{\partial \,C} {dA\,F}-\nabla \!  _VF\left( X \right)\,} \right|\le {\varepsilon  \over {\sup J(x)}}.   	\label {eq:e8}
\eeq

Let $\{ c_j, x_j\} _{j=1}^k$ be a $\delta$-fine RP partition of $[\,0,1]^n$. From the definition of the RP integral, Eq. \ref {eq:RP}, the approximating sum for the RP integral in Eq. \ref {eq:e9} is $\sum {\left| {c_j } \right|}J\left( {x_j} \right) V\!\left( {X_j} \right) \nabla \!  _VF\!\left( {X_j} \right)$. Then using Eq. \ref {eq:e8},
\beast
\lefteqn{\left| {\,\int_{\partial M} {dA\,F}\;-\;\sum\limits_j {\,\left| {\,c_j } \right|\,J\left( {x_j} \right)\;V(X_j)\nabla \!  _VF(X_j)\,}} \right| }  \\
&=&\left| {\,\sum\limits_j {\int_{\partial \,C_j} {dA\,F\;}}-\;\sum\limits_j {\,\left| {\,c_j} \right|\,J\left( {x_j} \right)\,V(X_j)\nabla \!  _VF(X_j)\,}} \right|  \\
&\le& \sum\limits_j {\left| {\,\int_{\partial \,C_j} {dA\,F\;}-\;\left| {\,c_j } \right|\,J\left( {x_j} \right)\,V(X_j)\nabla \!  _VF(X_j)\,} \right|}  \\
&=&\sum\limits_j {\left| {\,{{V^\dagger\!  (X_j)} \over {\left| {\,c_j} \right|\,J\left( {x_j} \right)}}\;\int_{\partial \,C_j} {dA\,F}\;-\;\nabla \!  _VF(X_j)\,} \right|\;}\,\left| {c_j} \right| J\left( {x_j} \right) \\
&\le& \sum\limits_j {{\varepsilon  \over {\sup J(x)}}\;}\left| {\,c_j} \right|\,J(x_j)\le \varepsilon .
\eeast
\vspace{.07in}

\noindent
{\bf Corollary}. (Cauchy-Goursat Theorem). \emph {Let $F$ be a multivector valued function defined on a $C^1$ $n$--chain $M$ in ${\mathbf R}^m$. Suppose that $F$ is continuous on $M$ and that $\nabla \!  _VF = 0$ on $M - \partial M$. Then $\int_{\partial M} {dA\,F=0}$ } 

\vspace{.07in}
(I elaborate a bit on the relationship to complex analysis. Suppose that the vector field ${\mathbf f} = ue_1 + ve_2$ is differentiable on an open set in ${\mathbf R}^2$. Define the complex function $f = v + iu$. Then by the Cauchy-Riemann equations,
\beast
\nabla {\mathbf f} =0 &\Leftrightarrow& \nabla \!  \cdot {\mathbf f} = 0\ \;\rm{and}\; \nabla \!  \wedge {\mathbf f}=0  \\
&\Leftrightarrow&u_x=-v_y \;\rm{and}\; u_y=v_x    \\
&\Leftrightarrow& f \;\rm{is }\;\rm{ analytic.}     
\eeast
Generalizing, a differentiable multivector field $F$ on an open set $U \subseteq M$ is called \emph{analytic}, or \emph{monogenic}, if $\nabla \!  _VF = 0$ on $U$. See \cite[p. 46]{bds}.)

\vspace{.15in}
{\bf V. Existence of }\mbox{\boldmath $\nabla _V \! F \ $}{\bf and its coordinate representation.}  
To show that $\nabla \!  _VF$ exists, we must show that the limit defining $\nabla\!  _VF$ in Eq. \ref {eq:e7} exists and is invariant under coordinate changes.

\vspace{.10in}
\noindent
{\bf Theorem.} \emph {
Let $M \subseteq {\mathbf R}^m$ be a $C^2$ $n$--chain and $F$ be a differentiable multi\-vector valued function defined on $M$. Then $\nabla \!  _VF$ exists. Moreover, using notation to be defined presently, }
\beq
\nabla \!  _VF(\varphi (x_0))=\sum\limits_{i=1}^n {{{V^\dagger\! (\varphi (x_0))\;\varphi '(x_0,a_i)} \over {J(x_0)}}\partial _iF(\varphi (x_0)).}  \label {eq:e10} 
\eeq
Here $\varphi$ is a singular $C^2$ $n$-cube. Let $h_i$ be the oriented hyperplane through $x_0$ with normal vector $+e_i$. Let $a_i$ be the tangent to $h_i$. For fixed $x$, the linear transformation $\varphi ' (x, \mathbf{\cdot})$ extends to an outermorphism from the tangent geometric algebra to ${\mathbf R}^n$ at $x$ to the tangent geometric algebra to $M$ at $\varphi (x)$ \cite[p. 165]{gc2}.  This defines $\varphi'(x_0, a_i)$.

Since $a_i$ is the unit tangent to $h_i$, $\varphi '(x_0, a_i)$ is a tangent to the surface $\varphi(h_i)$ at $\varphi(x_0)$ \cite[p. 165]{gc2}.  Thus, using $N = V{^\dagger\! }A$, the coefficient $N_i=V^\dagger\!  (\varphi (x_0))\;\varphi ' (x_0,a_i)/J(x_0)$
in Eq. (10) is a vector normal to $\varphi(h_i)$ at $\varphi(x_0)$. The $N_i$ need be neither orthogonal nor normalized. With this notation Eq. \ref {eq:e10} becomes $\nabla \!  _VF=\sum {N_i\partial _iF}$.

{\bf Proof.} For notational simplicity take $\varphi:[-1,1]^n \rightarrow M$ and $x_0 = 0$. Let $c \subseteq (-1,1)^n$ be a cube with $0 \in c$, of width $\epsilon$, and with sides $s_i^{\pm}$: on $s_i^+,\; x_i=\varepsilon _i^+>0$, and on $s_i^-, \;x_i=-\varepsilon _i^-<0$. Then $\varepsilon _i^++\varepsilon _i^-=\varepsilon $. And $\pm a_i$ is the tangent to $s_i^\pm$.

From the change of variables formula $dA = \varphi ' (x,da) = \left| da \right| \varphi'(x, a)$ \cite[p. 266]{gc2} and Eq. \ref {eq:e7},
\bea
\lefteqn {\nabla \! _VF\left( {\varphi  (0)} \right)={{V^\dagger\! \left( {\varphi (0)} \right)} \over {J\left(0\right)}} \; \times }  \label {eq:e11} \\
& & \lim \limits_{\varepsilon \to 0}\sum\limits_{i=1}^n {{1 \over {\varepsilon ^n}}\left[ {\int_{s_i^+} {\left| {\,da\,} \right|\;\varphi '  (x,a_i)\,F}(\varphi  (x))+\int_{s_i^-} {\left| {\,da\,} \right|\;\varphi ' (x,-a_i)\,\,F}(\varphi (x))} \right].}\, \nonumber
\eea

Since $F\circ \varphi $ is differentiable,
\beq
F(\varphi (x))=F(\varphi (0))+\sum\limits_{j=1}^n {\partial _jF(\varphi (0))}\,x_j+R(x), \label {eq:e12}
\eeq
where $\left| {R(x)} \right|/\left| x \right|\to 0$ as $\left| x \right | \to 0$.

We shall substitute separately the three terms on the right side of Eq. \ref {eq:e12} into the right side of Eq. \ref {eq:e11}. We shall find that the value of the first and third terms is 0 and that of the second term is the right side of Eq. \ref {eq:e10}, thus proving Eq. \ref {eq:e10}.

\textsc{First term.} Since $a_i = vn = ve_i$,
\bma
a_i=(-1)^{n-i}e_1\wedge \ldots \wedge \check{e}_i \wedge \ldots \wedge e_n, 
\ema
where $\check{e}_i$ indicates that $e_i$ is missing from the product.
Since the matrix element $\left[ {\varphi '} \right]_{ij}=\partial _j\varphi _i$, $\varphi ' \left( {x,e_j} \right)=\partial _j\varphi \left( x \right)$. Thus
\bma
\varphi '(x,a_i)=(-1)^{n-i}\partial _1\varphi (x)\wedge \ldots \wedge \partial _i \check{\varphi}(x)\wedge \ldots \wedge \partial _n\varphi (x).
\ema
Then (we justify setting the double sum to zero below)
\bea
\lefteqn {\sum\limits_{i=1}^n {\partial _i\varphi ' (x,a_i)} } \nonumber \\
& = &\sum\limits_{i=1}^n {\sum\limits_{\scriptstyle { {j=1} \atop {j\ne i} } } ^n  {(-1)^{n-i}\partial _1\varphi (x)\wedge \ldots \wedge \partial 
_{ji}\varphi (x)\wedge \ldots \wedge \partial _i \check{\varphi}(x)\wedge \ldots \wedge \partial _n\varphi (x)}}  \nonumber \label {eq:e13} \\
& = &0.
\eea
The terms in the double sum are paired by exchanging the $i$ and $j$ indices. The indicated term is paired with
\bma
(-1)^{n-j}\partial _1\varphi (x)\wedge \ldots \wedge \partial _j\check{\varphi} (x) \wedge \ldots \wedge \partial _{ij}\varphi (x)\wedge \ldots \wedge \partial _n\varphi (x).
\ema
Since $\varphi \in C^2$, $\partial _{ij}\varphi=\partial _{ji}\varphi$, and the paired terms differ by a factor of 
\bma
\left( {-1} \right)^{i-j-1}\left( {-1} \right)^{n-j}/\left( {-1} \right)^{n-i}=-1.
\ema
Thus the paired terms cancel and the double sum is zero. 

Let $x \in s_i^+$. Then $x_i = \epsilon _i^+$. Let $m_ix$ be the corresponding point on the opposite side $s_i^-$: $\left( m_ix \right)_i = -\epsilon _i^-$, and for $j \ne i$, $\left( m_ix \right)_j = x_j$. By the mean value theorem and our hypothesis $\varphi \in C^2$,
\bea
{{\varphi ' (x,a_i)\,-\varphi ' (m_ix,a_i)} \over \varepsilon }\,=\partial _i\varphi ' (x^*,a_i)=\partial _i\varphi ' (0,a_i)+S_i(x^*), \label {eq:e14}
\eea
where $x^*$ is between $x$ and $m_ix$ (and by abuse of notation $x^*$ is in general different for different components of the equation), and $\left| {S_i\left( {x^*} \right)} \right|\to 0$ as $\varepsilon \to 0$.

Substitute the first term on the right side of Eq. \ref {eq:e12} into Eq. \ref {eq:e11}. Notice that $\int_{s_i^+} {\left| {\,da\,} \right|}=\varepsilon ^{n-1}$. From Eqs. \ref {eq:e14} and \ref {eq:e13}, the sum in Eq. \ref {eq:e11} becomes
\beast
\lefteqn {\sum\limits_{i=1}^n {{1 \over {\varepsilon ^n}}\left[ {\int_{s_i^+} {\left| {\,da\,} \right|\;\varphi ' (x,a_i)\,F(\varphi (0))}+\int_{s_i^-} {\left| {\,da\,} \right|\;\varphi ' (x,-a_i)}\,F(\varphi (0))} \right]} } \\
& = &\sum\limits_{i=1}^n {{1 \over {\varepsilon ^n}}\int_{s_i^+} {\left| {\,da\,} \right|{{\varphi '(x,a_i)\,-\varphi ' (m_ix,a_i)} \over \varepsilon }\;\varepsilon }}\,F(\varphi (0))  \\
& \to & \sum\limits_{i=1}^n {\partial _i\varphi ' (0,a_i)}\,F(\varphi (0))=0.
\eeast

\textsc{Second term.} Substitute the second term on the right side of Eq. \ref {eq:e12} into Eq. \ref {eq:e11}. There are two cases in the resulting double sum: $i \ne j$ and $i = j$.

$i \ne j$. Recall that $\left(m_ix \right) _j=x_j$. Then
\pagebreak
\beast
\lefteqn { \left| {{1 \over {\,\varepsilon ^n}}} \left[ {\int_{s_i^+} {\left| {\,da\,} \right|\;\varphi '(x,a_i)\,\partial _iF(\varphi (0))\,x_j}} \; + \right. \right. }  \\
&  & \hspace{1in} \left. \left.{\int_{s_i^-} {\left| {\,da\,} \right|\;\varphi ' (x,-a_i)\,\partial _iF(\varphi (0))\,x_j} } \right] \right| \\
& = & \left| {\,{1 \over {\,\varepsilon ^n}}\int_{s_i^+} {\left| {\,da\,} \right|\,\left[ {\varphi ' (x,a_i)\,-\varphi ' (m_ix,a_i)} \right]}\;\partial _iF(\varphi (0))\,x_j\,} \right|  \\
& \le & \mathop {\sup }\limits_{x \in  s_i^+} \left| {\,\varphi ' (x,a_i)\,-\varphi ' (m_ix,a_i)\,} \right|\;\left| {\,\partial _iF(\varphi (0))\,} \right|\;\,\to \;0\,.
\eeast

$i = j$. Recall that on $s_i^\pm$, $x_i = \pm \varepsilon_i^\pm$. Then
\beast
\lefteqn{ {1 \over {\,\varepsilon ^n}} \left[{\int_{s_i^+} {\left| {\,da\,} \right|\;\varphi '(x,a_i)\,\partial _iF(\varphi (0))\,\varepsilon _i^+} \; +} \right. } \\
& & \hspace{1in} \left. {\int_{s_i^-} {\left| {\,da\,} \right|\;\varphi '(x,-a_i)\,\partial _iF(\varphi (0))\,(-\varepsilon _i^-)}} \right] \\
& = &{1 \over {\varepsilon ^n}}\int_{s_i^+} {\left| {\,da\,} \right| \; \left[ {\varphi '(x,a_i)\,\varepsilon _i^++\varphi '(m_ix,a_i)\,\varepsilon _i^-} \right]}\;\partial _iF(\varphi (0)) \\
& = & \frac{1}{\varepsilon ^n} \int_{s_i^+} {\left| {\,da\,} \right| \, \left\{  {\,\left[ {\varphi '(x,a_i)-\varphi '(0,a_i)} \right]\,\varepsilon _i^+ } \right.  } \; + \\
& & \hspace{1in} { \left. { \left[ {\varphi '(m_ix,a_i)-\varphi '(0,a_i)} \right]\,\varepsilon _i^-+\varphi '(0,a_i)\,\,\varepsilon \,} \right\}\;}\partial _iF(\varphi (0)) \\
& \to & \varphi '(0,a_i)\;\partial _iF(\varphi (0))\,,
\eeast
which when multiplied by $V^\dagger\! \left( {\varphi \left( 0 \right)} \right)/J\left( 0 \right)$ from Eq. \ref {eq:e11} is the $i^{\rm th}$ summand in Eq. \ref {eq:e10}.

\textsc{Third term.} Substitute $R(x)$ into the first (for example) term in a summand on the right side of Eq. \ref {eq:e11}. Then since $ {\left| R(x) \right|} / {\left| x \right| } \to 0$ as $\left| x \right | \to 0$,
\bma
{1 \over {\varepsilon ^n}} 
\left| \int_{s_i^+} {\left| da \right|        \varphi'(x,a_i)\,R(x) \, }\right|
\le {1 \over {\varepsilon ^n}}
       \int_{s_i^+} {\left| da \right| \left| \varphi'(x,a_i)\right| \sqrt{n} \varepsilon \;
         {{\left| R(x) \right|} \over {\left| x \right|}  }} \to \;0\,.
\ema 
(From \cite[pp. 13-14]{gc2} we see that $|x + y| \le |x| + |y|$ and that if $x$ is simple, then $|xy| = |x| \, |y|$. This justifies the inequality above.)

\vspace{.1in}
The coordinate representation for $\nabla F$, Eq. \ref {eq:e10}, is now proved. Our remaining task is to show that the representation is invariant under a coordinate change. Suppose then that $\psi :\left[ {-1,1} \right]^n\to M$ is a second $C^2$ singular {\it{n}}--cube with, again for notational convenience, $\psi(0) = \varphi(0)$. For this {\it{n}}--cube use the notations $y, f_i, b_i,$ and $w$ corresponding to $x, e_i, a_i,$ and $v$ above. Define $g(y) = x$ by $g=\varphi ^{-1}\circ \psi $.

Reversing $v^\dagger a_j = e_j$ gives $a_j^\dagger v = e_j$, i.e., $a_jv = 
(-1)^{(n-1)n/2}e_i$. Similarly, $b_iw = (-1)^{(n-1)n/2}f_i$. Let $\bar g$ be the adjoint map to $g'$. The matrix element $[\bar g]_{ij}=[g]_{ji}={{\partial x_j} \mathord{\left/ {\vphantom {{\partial x_j} {\partial y_i}}} \right. } {\partial y_i}}$. Then
\bea
 \bar g(0,\;a_j\,v) &=& (-1)^{(n-1)n/2}\,\bar g(0,\;e_j)  \nonumber \\
&=&(-1)^{(n-1)n/2}\sum\limits_{i=1}^n {{{\partial \,x_j} \over {\partial \,y_i}}}\,f_i=\sum\limits_{i=1}^n {{{\partial \,x_j} \over {\partial \,y_i}}}\,b_i\,w.  \label {eq:e15}
\eea
The outermorphisms $g'$ and $\bar {g}$ are related by a duality \cite[Eq. 3-120b]{gc2}. By this duality the equality of the outermost members of Eq. \ref {eq:e15} is equivalent to
\bma
{1 \over {J_g(0)}}\:g'\!\left( {0,\;\sum\limits_{i=1}^n {{{\partial \,x_j} \over {\partial \,y_i}}}\,b_i} \right)=a_j.
\ema

With this identity we show that the right side of Eq. \ref {eq:e10} is invariant under a coordinate change:
\beast
\lefteqn{\sum\limits_{i=1}^n{{{V^\dagger\!  (\psi (0))\;\psi '(0,b_i)} \over {J_\psi (0)}}\;\partial _{y_i}F(\psi (0))} }  \\
&=& \sum\limits_{i=1}^n {\left\{ {{{V^\dagger\!  (\varphi (0))\;(\varphi \circ g)'(0,b_i)} \over {J_{\varphi \circ g}(0)}}\; 
\sum\limits_{j=1}^n {\left[ {\partial _{x_j}F(\varphi (0)){{\partial \,x_j} \over {\partial \,y_i}}} \right]}} \right\}}  \\
&=& \sum\limits_{j=1}^n {\left\{ {{{V^\dagger\!  (\varphi (0))} \over {J_\varphi  (0)}}\;\varphi' \! \left[ {0,\;{1 \over {J_g(0)}}\:g'\!\left( {0,\;\sum\limits_{i=1}^n {{{\partial \,x_j} \over {\partial \,y_i}}b_i}} \right)} \right]\;\partial _{x_j}F(\varphi (0))} \right\}} \\
&=&  \sum\limits_{j=1}^n {{{V^\dagger\!  (\varphi (0))\;\varphi ' (0,a_j)} \over {J_\varphi (0)}}\partial _{x_j}F(\varphi (0)).}  \\
\eeast

\pagebreak

\end {document}